\documentclass[b5paper]{article}
\usepackage{geometry}
\usepackage{amsthm,amssymb,amsmath}
\usepackage{url}

\begin{document}

\title{A note on ``MLE in logistic regression with a diverging dimension''}
\author{Huiming Zhang}
\date{\footnotesize{\today}}
\maketitle

\begin{abstract}
This short note is to point the reader to notice that the proof of high dimensional
asymptotic normality of MLE estimator for logistic regression under the regime $p_n=o(n)$ given in paper: ``Maximum likelihood estimation in logistic regression models with a diverging number of covariates. \emph{Electronic Journal of Statistics}, 6, 1838-1846." is wrong.
\\
\\
\textbf{Keyword}: high dimensional logistic regression; generalized linear models; asymptotic normality.
\end{abstract}
~
\\
\\
In order to maintain the preciseness of the statistical scientific record, I write and post this notes on arxiv, as an intention for avoiding misleading the reader.

Under mild conditions, it seems that \cite{Liang} gives a concise proof that
the MLE for logistic regression is asymptotically normality when the number of covariates $p$ goes to infinity with the sample size $n$ satisfying $p_n=o(n)$. \cite{Liang} claimed that their results sharpen the existing results of asymptotic normality, for example, \cite{Portnoy} studied the MLE estimator for GLM with dimension rate increasing $p_n = o({n^{2/3} })$ and \cite{Wang} analysed the GEE estimator for logistic regression as $p_n = o({n^{1/3} })$. The result of \cite{Liang} was cited by some recent papers, such as \cite{Sur} for the high-dimensional likelihood ratio test in logistic regression when $\frac{p_n}{n} < \frac{1}{2}$.

Nevertheless, it can be carefully seen that Lemma 3 in \cite{Liang}, which is adapted from claims (18) and (19) in \cite{Yin} for quasi-likelihood estimates, is not true.  We restate Lemma 3 and its first part of the proof in \cite{Liang} as below:
\begin{quotation}
\footnotesize{\textbf{Lemma 3} (\cite{Liang}). Under the conditions of Theorem 1, we have
\[\mathop {\sup }\limits_{\beta  \in {N_n}(\delta )} \left| {{u^T}G_n^{ - 1/2}({\beta _0}){Q_n}(\beta )G_n^{ - 1/2}({\beta _0})u - 1} \right| \to 0.~~~~({u^T}u = 1)\]
where ${Q_n}(\beta ): = \frac{{\partial {L_n}(\beta )}}{{\partial {\beta ^T}}}$, ${G_n}({\beta _0}) = \sum\limits_{i = 1}^n {{x_i}} H(x_i^T{\beta _0})x_i^T$ and ${N_n}(\delta ): = {\rm{\{ }}\beta : \parallel G_n^{ - 1/2}({\beta _0})(\beta  - {\beta _0})\parallel  \le \delta {\rm{\} }}$.
\begin{proof}
Let ${\varepsilon _i}: = {y_i} - h(x_i^T\beta _0),\sigma _i^2 = {\mathop{\rm var}} {\varepsilon _i}$. A direct calculation yields
\begin{equation} \label{eq:dc}
{u^T}G_n^{ - 1/2}({\beta _0}){Q_n}(\beta )G_n^{ - 1/2}({\beta _0})u - 1{\rm{ = }}{A_n}(\beta ) - {B_n}({\beta _0}) - {C_n}(\beta )
\end{equation}
where ${A_n}(\beta ) = G_n^{ - 1/2}({\beta _0}){G_n}(\beta )G_n^{ - 1/2}({\beta _0}) - 1$, ${B_n}({\beta _0}) = \sum\limits_{i = 1}^n {{u^T}G_n^{ - 1/2}({\beta _0}){x_i}} x_i^TG_n^{ - 1/2}({\beta _0}){u^T}{\varepsilon _i}$, \\ ${C_n}(\beta ) = \sum\limits_{i = 1}^n {{u^T}G_n^{ - 1/2}({\beta _0}){x_i}} x_i^TG_n^{ - 1/2}({\beta _0}){u^T}[h(x_i^T{\beta _0}) - h(x_i^T\beta )].\cdots$.
\end{proof}}
\end{quotation}
The so-called ``direct calculation" of verifying the decomposition (1) (modified from decomposition (21) in \cite{Yin}) is not true. Notice that ${G_n}(\beta )$ and ${Q_n}(\beta )$ is frankly equal for the setting in \cite{Liang} (first paragraph of p1840), thus we must have
\[{B_n}({\beta _0}) + {C_n}(\beta ) = \sum\limits_{i = 1}^n {{u^T}G_n^{ - 1/2}({\beta _0}){x_i}} x_i^TG_n^{ - 1/2}({\beta _0}){u^T}[{y_i} - h(x_i^T\beta )] \equiv 0,\]
which is impossible.

In fact, the definition ${Q_n}(\beta ): = \frac{{\partial {L_n}(\beta )}}{{\partial {\beta ^T}}}$ in Lemma 3 of \cite{Liang} is straightly borrowed from (21) in \cite{Yin} for partial derivatives of log quasi-likelihood w.r.t. $\beta$
\[L_n^{quasi}(\beta ): = \sum\limits_{i = 1}^n {{x_i}} H(x_i^T\beta ){[\Sigma (x_i^T\beta )]^{ - 1}}[{y_i} - h(x_i^T\beta )]\]
where $H(t):= \frac{{dh(t)}}{{dt}},\Sigma (x_i^T\beta ): = Co{v_\beta }({y_i})$.

Directly applying decomposition (21) in \cite{Yin} is not make sense, since the partial derivatives of log quasi-likelihood w.r.t. $\beta $ depend on the responses $\{ {y_i}\} _{i = 1}^n$, and the responses in partial derivatives of log likelihood
\[{L_n^{mle}}(\beta ): = \sum\limits_{i = 1}^n {{x_i}} [{y_i} - h(x_i^T\beta )]\]
are cancelled, i.e. ${Q_n}(\beta ) = \sum\limits_{i = 1}^n {{x_i}} H(x_i^T\beta )x_i^T$ which is not random in Lemma 3 of \cite{Liang}, since the covariates are assumed to be deterministic. Actually the symbol ``$Q$"  means ``quasi-". However, \cite{Liang} confusedly utilized it for log likelihood of logistic regression.

For rest of the proof, the technique of asymptotic derivation is almost the same as \cite{Yin} whose matrix computations and mathematical arguments (such as the local inverse function theorem) play an essential role.

In the last paragraph of \cite{Liang}, they say: ``We believe that the procedure can be extended to other generalized linear models and similar theoretical results may be established with straightforward derivations. One potential complication for other generalized linear models is that the response $y$ may not be bounded as in logistic regression models. Other possible extensions are to the Cox model, robust regression, and procedures based on quasi-likelihood functions. Further effort is needed to build up similar procedure and theoretical results under these settings."

6 years after publication from Google Scholar citation, there were not any related papers which extended their results $p_n=o(n)$ to other generalized linear models with straightforward derivations.

For example, The Theorem 2 of \cite{Khaplanov} who cited \cite{Liang}. This Russian paper extend \cite{Wang} to multivariate logistic regression with a diverging number of covariates, it just obtains asymptotic normality of MLE with $p_n = o({n^{1/3} })$. The proof techniques in \cite{Khaplanov} are also borrowed from \cite{Yin}.

\bigskip

\noindent\textit{School of Mathematical Sciences,
Peking University, Beijing, P.R.China\\
zhanghuiming@pku.edu.cn}

\bigskip

\end{document}